\documentclass{article}
\usepackage{spconf,amsmath,epsfig}
\RequirePackage{color}
\usepackage{amssymb}
\usepackage{amsthm}
\usepackage{amsfonts}
\usepackage{bbm}
\usepackage{mathrsfs}

\def\Re{\mathbb{R}}
\def\ones{\mathbbm{1}}
\newcommand{\ints}{{\mathbb{Z}}}
\def \Poisson{{{\rm Poisson}}}
\def\deq{\triangleq}
\def\hf{\widehat{f}}
\def\tf{\widetilde{f}}

\def\argmin{\mathop{\arg \min}}
\def\pen{{\rm pen}}
\def\sP{\mathcal P}

\newcommand{\squishlist}{
   \begin{list}{$\bullet$}
    { \setlength{\itemsep}{0pt}      \setlength{\parsep}{0pt}
      \setlength{\topsep}{0pt}       \setlength{\partopsep}{0pt}
      \setlength{\leftmargin}{1.5em} \setlength{\labelwidth}{1em}
      \setlength{\labelsep}{0.5em} } }
\newcommand{\squishend}{
    \end{list}  }
\usepackage{booktabs}

\title{Sparse Poisson Intensity Reconstruction Algorithms}
\name{Zachary T. Harmany, Roummel F. Marcia, and Rebecca M. Willett
  \thanks{This work was supported by NSF CAREER Award No.\
    CCF-06-43947, NSF Award No.\ DMS-08-11062, DARPA Grant
    No. HR0011-07-1-003, and AFRL Gran No. FA8650-07-D-1221.}}
\address{Duke University, Department of Electrical and Computer
  Engineering, Durham, NC 27708 USA}
\begin{document}
\ninept
\maketitle
\begin{abstract}
  \small \sloppypar{ The observations in many applications consist of
    counts of discrete events, such as photons hitting a dector, which
    cannot be effectively modeled using an additive bounded or
    Gaussian noise model, and instead require a Poisson noise
    model. As a result, accurate reconstruction of a spatially or
    temporally distributed phenomenon ($f$) from Poisson data ($y$)
    cannot be accomplished by minimizing a conventional
    $\ell_2-\ell_1$ objective function. The problem addressed in this
    paper is the estimation of $f$ from $y$ in an inverse problem
    setting, where (a) the number of unknowns may potentially be
    larger than the number of observations and (b) $f$ admits a sparse
    approximation in some basis.  The optimization formulation
    considered in this paper uses a negative Poisson log-likelihood
    objective function with nonnegativity constraints (since Poisson
    intensities are naturally nonnegative).  This paper describes
    computational methods for solving the constrained sparse Poisson
    inverse problem. In particular, the proposed approach incorporates
    key ideas of using quadratic separable approximations to the
    objective function at each iteration and computationally efficient
    partition-based multiscale estimation methods. }
\end{abstract}

\section{Introduction}
\label{sec:intro}

In a variety of applications, ranging from nuclear medicine to night
vision and from astronomy to traffic analysis, data are collected by
counting a series of discrete events, such as photons hitting a
detector or vehicles passing a sensor. The measurements are often
inherently noisy due to low count levels, and we wish to reconstruct
salient features of the underlying phenomenon from these noisy
measurements as accurately as possible. The inhomogeneous Poisson
process model \cite{snyder} has been used effectively in many such
contexts. Under the Poisson assumption, we can write our observation
model as
\begin{equation}
y \sim \Poisson(Af), \label{eq:Poisson}
\end{equation}
where $f \in \Re_+^{m}$ is the signal or image of interest, $A \in
\Re_+^{N \times m}$ linearly projects the scene onto a $N$-dimensional
set of observations, and $y \in \ints_+^N$ is a length-$N$ vector of
observed Poisson counts. 

The problem addressed in this paper is the estimation of $f$ from $y$
in a compressed sensing context, when (a) the number of unknowns, $m$,
may be  larger than the number of observations, $N$, and (b) $f$ is
sparse or compressible in some basis $W$ (i.e. $f = W\theta$ and
$\theta$ is sparse). This challenging problem has clear connections to
``compressed sensing'' (CS) (e.g., \cite{RIP}), but arises in a number of other settings as
well, such as tomographic reconstruction in nuclear medicine,
superresolution image reconstruction in astronomy, and deblurring in
confocal microscopy.

In recent work \cite{WilRag09}, we explored some of the theoretical
challenges associated with CS in a Poisson noise
setting, and in particular highlighted two key differences between the
conventional CS problem and the Poisson CS problem:
\begin{enumerate}
\item unlike many sensing matrices in the CS literature, the matrix
  $A$ must contain all nonnegative elements, and
\item the intensity $f$, and any estimate of $f$, must be nonnegative.
\end{enumerate}
Not only are these differences significant from a theoretical
standpoint, but they pose computational challenges as
well. 
In this paper, we consider algorithms for computing as estimate $\hf$
of $f$ from observations $y$ according to
\begin{eqnarray*}
	\hf &\deq& \argmin_{\tf \in \Gamma} -\log \ p \! \left(y \left| A\tf 
          \right. \right) + \pen (\tf \, ) \label{eq:objective}\\
	&=& \argmin_{\tf \in \Gamma} \sum_{i=1}^N  \left \{ 
		\left(A\tf \, \right)_i -  y_i \log
		\left[\left(A\tf\,\right)_i\right] 
		\right \}
	+ \pen (\tf \, ) ,\nonumber 
\end{eqnarray*}
where 
\squishlist
\item $p \! \left(y \left| A\tf \right.\right)$ is the Poisson likelihood of $y$ given intensity
$A\tf$,
\item $\left(A\tf\,\right)_i$ is the $i^{\rm th}$ element of the
vector $A\tf$,
\item $\Gamma$ is a collection of nonnegative estimators such that
  $\tf \ge 0$ for all $\tf \in \Gamma$, and
\item $\pen ( \tf \, )$ is a penalty term which measures the
  sparsity of $\tf$ in some basis. 
\squishend
This paper explores computational methods for solving
(\ref{eq:objective}).  The nonnegativity of $f$ and $A$ results in
challenging optimization problems. In particular, the restriction that
$f$ is nonnegative introduces a set of inequality constraints into the
minimization setup; as shown in \cite{fesslerPoissonCS}, these
constraints are simple to satisfy when $f$ is sparse in the canonical
basis, but they introduce significant challenges when enforcing
sparsity in an arbitrary basis.  We will consider several variants of
the penaly term, including $\|f\|_1$, $\|W^T\! \! f\|_1$ for some arbitrary
orthonormal basis $W$, and a complexity regularization penalty based
upon recursive dyadic partitions. We refer to our approach as SPIRAL
(Sparse Poisson Intensity Reconstruction ALgorithm).

\sloppypar{
The paper is organized as follows.  We first consider related optimization 
methods for solving the compressed sensing problem in Sec.\ \ref{sec:regular}. 
Next, in Sec.\ \ref{sec:penalties}, we consider three methods for
solving  the constrained penalized Poisson likelihood estimation
problem with different penalty terms.  The first two approaches use an $\ell_1$
penalty while the third uses partition-based penalty functions. 
Experimental results comparing the proposed approaches with 
expectation-maximization (EM) algorithms 
are presented in Sec.~\ref{sec:simulations}.  
}

\section{Sparsity-regularized optimization}
\label{sec:regular}

\sloppypar{
Without any nonnegativity constraints on $f$, 
the Poisson minimization problem (\ref{eq:objective}) 
can be solved using the SpaRSA algorithm of Figueireido et al.\ \cite{sparsa}, which 
solves a sequence of minimization problems using 
quadratic approximations to the log penalty function:
}
\begin{equation}	\label{eq:sparsa_f}
	f^{k+1} \deq
	\underset{f \in \Re^m}{\arg \min } \ \ 
	\frac{1}{2} \| f - s^k \|_2^2 + \frac{1}{\alpha_k} \pen (f),
\end{equation}
where 
\begin{equation}	\label{eq:sk}
	s^k = f^k - \frac{1}{\alpha_k} \nabla F(f^k),	
\end{equation}
$\alpha_k$ is a positive scalar chosen using Barzilai-Borwein (spectral) methods (see \cite{sparsa} for details), and
$$F(f) \deq -\log p(y|Af)$$
 is the negative log Poisson likelihood in
(\ref{eq:objective}).
If the penalty term in (\ref{eq:sparsa_f}) is separable, i.e., it is the sum of functions with each function
depending only on one component of $f$, then (\ref{eq:sparsa_f}) can be solved component-wise,
making the problem (\ref{eq:sparsa_f}) relatively easy to solve.

Much of existing CS optimization literature (e.g., \cite{l1-magic, GPSR, bregman})
focuses on penalty functions where $\pen(f) \propto \|W^T\!\!f\|_1$, where $W$ is some
orthonormal basis, such as a wavelet basis, and $\theta \deq W^T\!\!f$ are
the expansion coefficients in that basis.   Also, nonnegativity in the signal $f = W\theta \ge 0$
is not necessarily enforced.   
In this paper, we develop three approaches to solving the optimization problem in (\ref{eq:objective}).
In all three approaches, we require $f$ to be nonnegative.  
The first approach assumes that the image is sparse in the canonical basis, i.e., 
$W = I$, while the second involves the more general orthogonal matrix $W$.  The third uses
partition-based denoising methods, described in detail in Sec.~\ref{sec:partition}.

\section{nonnegative regularized least-squares subproblem}
\label{sec:penalties}

\subsection{Sparsity in the canonical basis} Let 
$\pen(f) = \tau \| f \|_1$, where $\tau > 0$ is a
regularization parameter.  The minimization problem
(\ref{eq:sparsa_f}) with this penalty term and that is subject to
nonnegativity constraints on $f$ has the following analytic solution:
\begin{equation*}
	f^{k+1} = \left [ s^k - \frac{\tau}{\alpha_k} \ones \right ]_+
\end{equation*}
where the operation $[ \enspace \cdot \enspace ]_+ =  \max \{ \  \cdot
\ , 0 \}$ is to be understood component-wise.
Thus solving (\ref{eq:objective}) subject to nonnegativity constraints with an $\ell_1$ penalty function and
with a solution that is sparse in the canonical basis is straightforward to solve.
An alternative algorithm for solving this Poisson inverse
problem with sparsity in the canonical basis was also explored in the
recent literature \cite{fesslerPoissonCS}.

\subsection{Sparsity in an arbitrary orthonormal basis}
\label{sec:basis_sparsity}
Now suppose that the signal of interest is sparse in some other basis.  Then
the $\ell_1$ penalty term is given by
$$
\pen(f) \deq \tau \|W^T\!\!f\|_1 =  \tau |\theta\|_1,
$$
where $\theta \deq Wf$ for some orthonormal basis $W$, and $\tau > 0$ is some scalar.
 When the reconstruction $\hat{f} = W\hat{\theta}$ must be nonnegative (i.e.,
$W\hat{\theta} \ge 0$), the minimization problem
\begin{align}	\nonumber
	\theta^{k+1} \deq \ &
	\underset{\theta \in \Re^m}{\arg \min } \quad
	\phi^{k}(\theta) \deq \frac{1}{2} \| \theta - s^k \|_2^2 + \frac{\tau}{\alpha_k} \| \theta \|_1, \\
	& \textrm{subject to \ } W\theta \ge 0
	\label{eq:sparsa_con} 
\end{align}
no longer has an analytic solution necessarily.  
We can solve this minimization problem by solving its dual.
First, we reformulate 
(\ref{eq:sparsa_con}) so that its objective function $\phi^k(\theta)$ is differentiable 
by defining $u, v \ge 0$ such that $\theta = u - v$.  The minimization problem (\ref{eq:sparsa_con})
becomes
\begin{align}	
	\nonumber
	(u^{k+1}, v^{k+1}) \deq \ &
	\underset{u,v \in \Re^m}{\arg \min } \quad
	\frac{1}{2} \| u -  v - s^k \|_2^2 + \frac{\tau}{\alpha_k} \ones^T\!(u+v) \\
	\label{eq:primal}
	& \textrm{subject to \ } u,v \ge 0, \quad W(u-v) \ge 0, \\[-.5cm] \nonumber
\end{align}
which has twice as many parameters and has additional nonnegativity constraints on the new parameters, 
but now has a differentiable objective function.
The Lagrange dual problem associated with this problem is given by
\vspace{-.2cm}
\begin{eqnarray}	\nonumber
	&& \underset{\lambda, \gamma \in \Re^m}{\textrm{minimize}} \ \ 
	h(\lambda, \gamma) \deq
	\frac12 \| s^k + \gamma  + W^T\! \lambda  \|_2^2 - \frac12 \| s^k \|_2^2  \quad \\ 
	&& \textrm{subject to \ } \lambda \ge 0, \ \ \  - \frac{\tau}{\alpha_k} \ones \le \gamma \le \frac{\tau}{\alpha_k} \ones 
	\label{eq:dual} \\[-.3cm] \nonumber
\end{eqnarray}
and at the optimal values $\gamma^{\star}$ and
$\lambda^{\star}$, the primal iterate $\theta^{k+1} \deq u^{k+1} - v^{k+1}$ is given by
$$
	\theta^{k+1} = s^k + \gamma^{\star} + W^T\! \lambda^{\star}.
$$
We note that the minimizers of the primal problem (\ref{eq:primal})
and its dual (\ref{eq:dual}) satisfy $\phi^k(\theta^{k+1}) = -h(\gamma^{\star}, \lambda^{\star})$
since (\ref{eq:primal}) satisfies (a weakened) Slater's condition \cite{convexoptimization}.
In addition, the function $-h(\gamma, \lambda)$ is a lower bound on $\phi^{k}(\theta)$
at any dual feasible point.

The objective function of (\ref{eq:dual}) can be minimized by alternatingly solving for $\lambda $ and $\gamma$,
which is accomplished by taking the partial derivatives of $h(\lambda, \gamma)$ and setting them to zero. 
Each component is then constrained to satisfy the bounds in (\ref{eq:dual}).  At the $j^{\rm th}$ iteration, 
the variables can, thus, be defined as follows:
\begin{eqnarray*}
	\gamma^{(j)} \!\!\! &=& \!\! \textrm{mid}  \left \{  -\frac{\tau}{\alpha_k} \ones, -s^k -W^T \! \lambda^{(j-1)},
			\frac{\tau}{\alpha_k} \ones \right \} \\
	\lambda^{(j)} \!\!\! &=& \!\! \left [ -W \left (s^k + \gamma^{(j)} \right ) \right]_+,
\end{eqnarray*}
where the operator mid($a,b,c$) chooses the middle value of the three arguments component-wise.  
Note that at the end of each iteration $j$, the 
approximate solution $\theta^{(j)} \deq s^k + \gamma^{(j)} + W^T \lambda^{(j)}$ to (\ref{eq:sparsa_con}) 
is feasible with respect to the constraint $W\theta \ge 0$:
\begin{eqnarray*}
	W\theta^{(j)} &=& 
	Ws^k + W\gamma^{(j)}  + \lambda^{(j)}  \\
	&=& \left [ W \left (s^k + \gamma^{(j)} \right )  \right ]_+ \ge 0.
\end{eqnarray*}

Unfortunately, there are several important disadvantages associated
with the above approach. First, the solution to the subproblem must
itself be computed by an iterative algorithm -- so while the problem
may be solvable, it must by solved at each iterate, and the resulting overall algorithm
may not be fast, particularly for large problems. Furthermore, since
we rely on an iterative solver in this subproblem, our computed
solution may not be the true optimal point, particularly if we use a
gentle convergence requirement to limit computation time.  Further
study is needed to understand the impact
of computing inaccurate solutions to this subproblem on the performance
of the overall algorithm.  

\subsection{Partition-based penalties} \label{sec:partition}

In the special case where the signal of interest is known to be smooth
or piecewise smooth in the canonical basis (i.e. is compressible in a
wavelet basis), we can formulate a penalty function which is a useful
alternative to the $\ell_1$ norm of the wavelet coefficient vector. In
particular, we can build on the framework of {\em recursive dyadic
  partitions (RDP)}, which we summarize here and are described in detail in
\cite{willett:density,willett:jsac04}. Let $\sP$ be the
class of all recursive dyadic partitions of $[0,1]$ where each
interval in the partition has length at least $1/m$, and let $P \in
\sP$ be a candidate partition. The intensity on $P$, denoted $f(P)$,
is calculated using a nonnegative least-squares method to fit a model
(such as a constant or polynomial) to $s^k$ in (\ref{eq:sk}) on each interval in the RDP.
Furthermore, a penalty can be assigned to the resulting
estimator which is proportional to $|P|$, the number of intervals in $P$.
Thus we set
\vspace{-.15cm}
\begin{align*}
\widehat{P} &= \argmin_{P \in \sP} \frac{1}{2} \|f(P) - s^k\|_2^2 + \tau |P|\\
\hf &= f(\widehat{P}).
\end{align*}
A search over $\sP$ can be computed quickly using a dynamic program. When using
constant partition interval models , the nonnegative
least-squares fits can be computed non-iteratively in each interval by
simply using the maximum of the average of $s^k$ in that interval and
zero. Because of this, enforcing the constraints is trivial and can be
accomplished very quickly.

It can be shown that partition-based denoising methods such as this
are closely related to Haar wavelet denoising with an important
hereditary constraint placed on the thresholded coefficients---if a
parent coefficient is thresholded, then its children coefficients must
also be thresholded \cite{willett:density}. This constraint is akin to
wavelet-tree ideas which exploit persistence of significant wavelet
coefficients across scales and have recently been shown highly useful
in compressed sensing settings \cite{baraniukModelBasedCS}.

As we mention in Sec.~\ref{sec:simulations}
below, using constant model fits makes it easy to satisfy
nonnegativity constraints and results in a very fast non-iterative
algorithm, but has the disadvantage of yielding piecewise constant
estimators. However, a cycle-spun translation-invariant version of this approach
can be implemented with high computational efficiency
\cite{willett:density} and be used for solving this nonnegative regularized
least-squares subproblem that results in a much smoother estimator.

\vspace{-.15cm}

\section{Numerical simulations}
\label{sec:simulations}

\begin{figure*}[t!]
\centering
\begin{tabular}{cc}  
	\includegraphics[width=7.3cm]{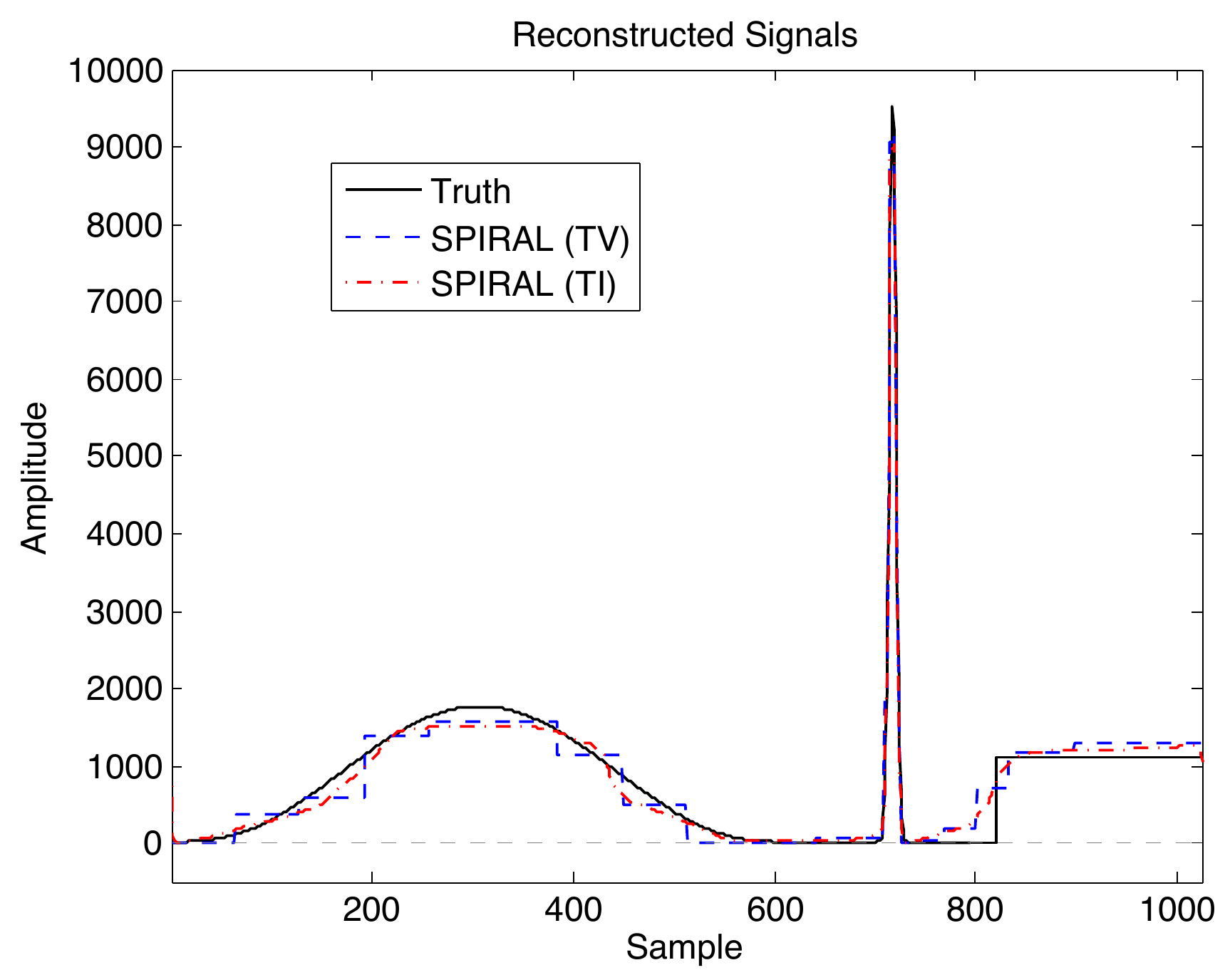} &
	\includegraphics[width=7.3cm]{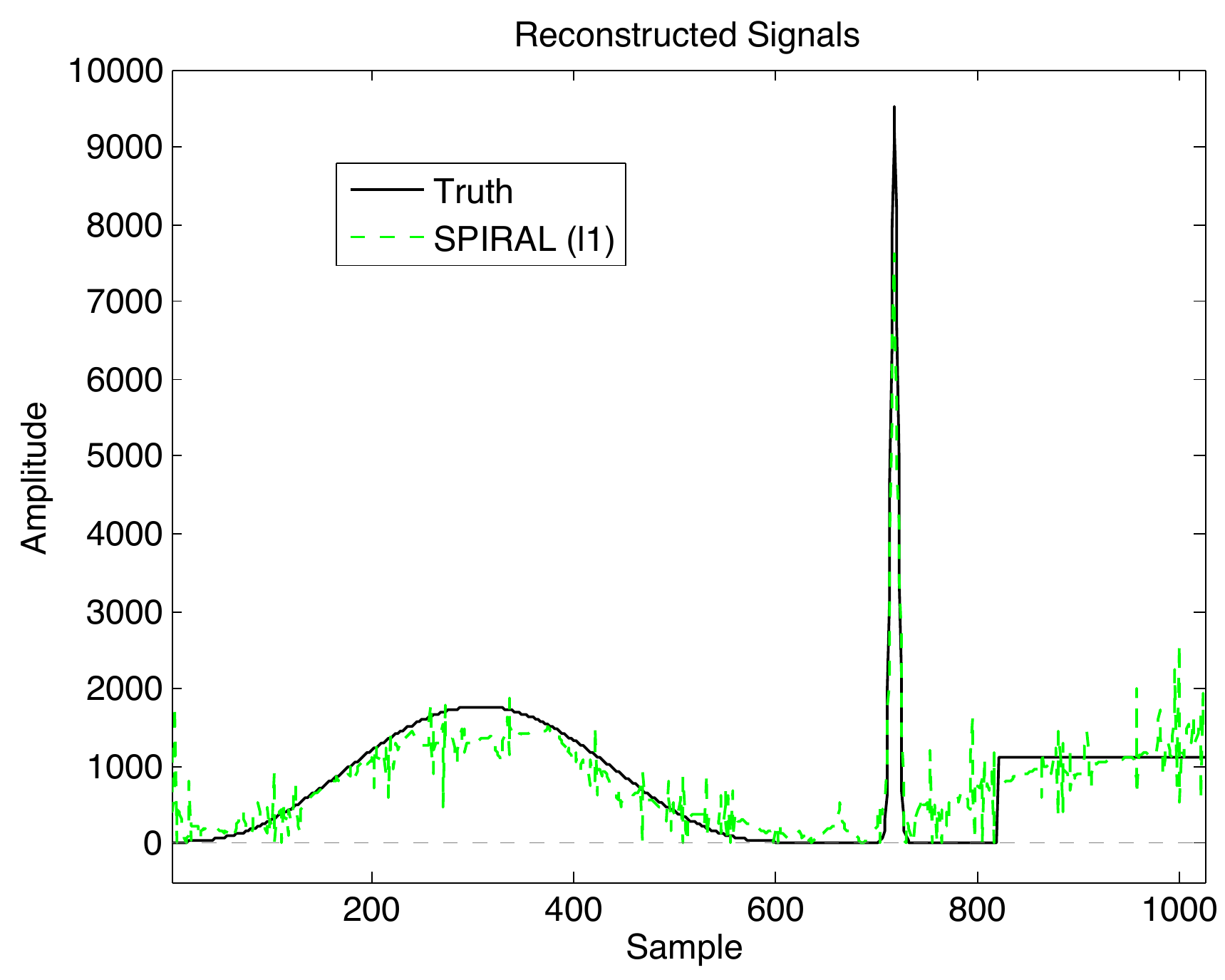} \\
  (a) & (b) \\
  \includegraphics[width=7.3cm]{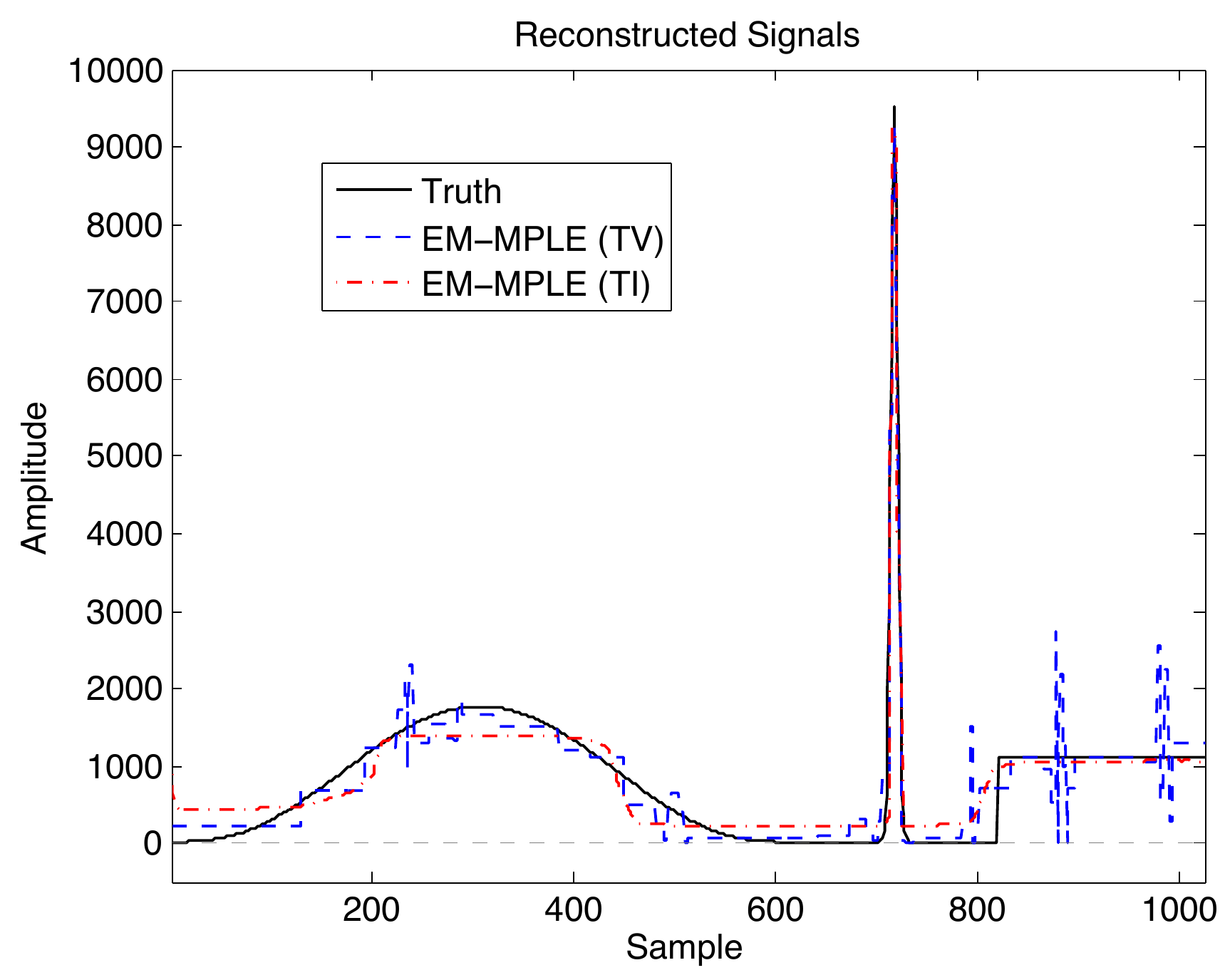} &
  \includegraphics[width=7.3cm]{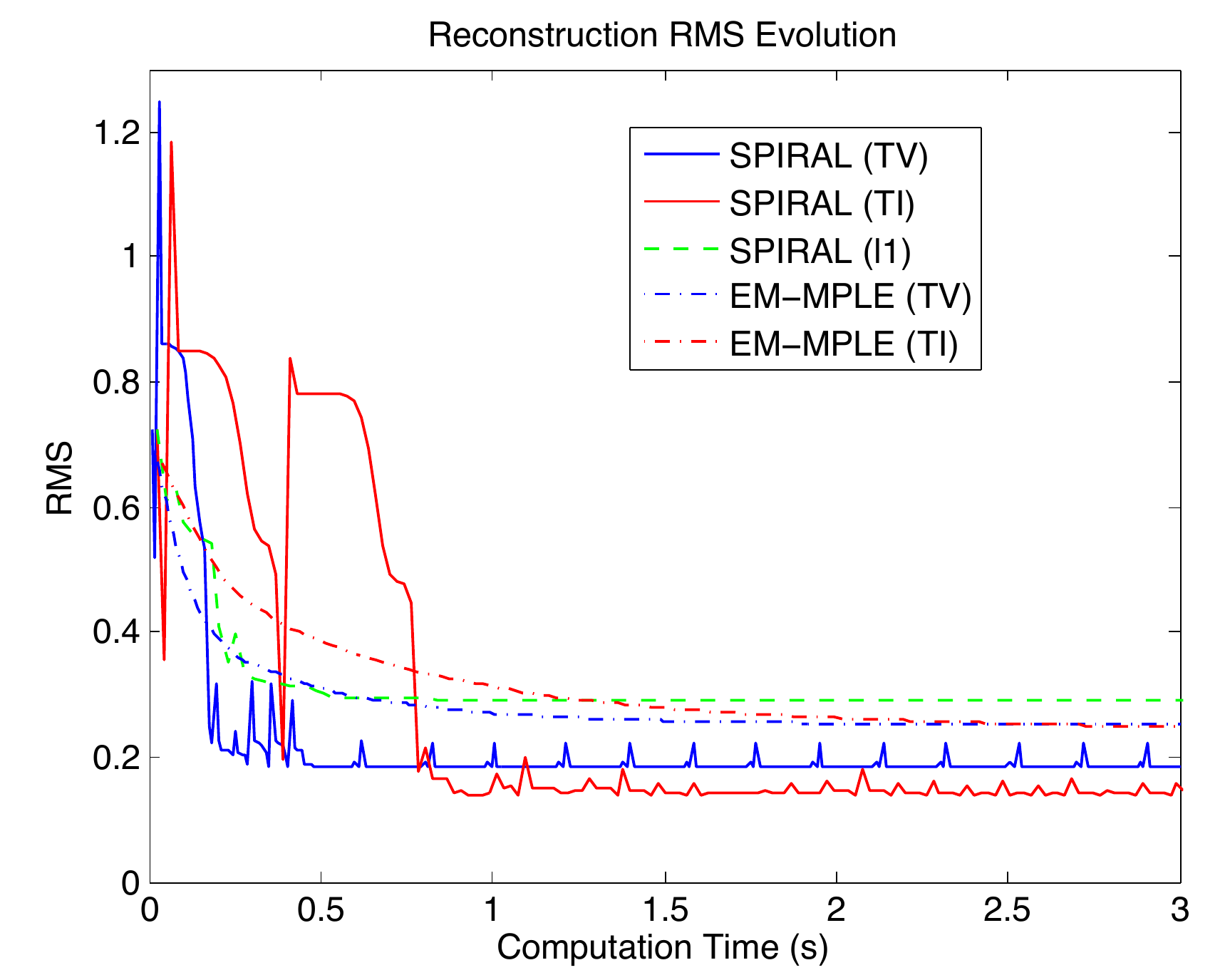} \\
  (c) & (d) 
  \end{tabular}
  \caption{Reconstruction results using(a) the SPIRAL algorithm with
  		 translationally variant (TV) and translationally invariant (TI) partitions,
		 (b) the SPIRAL algorithm with an $\ell_1$ penalty, and (c) the EM-MPLE algorithm with 
  		translationally variant (TV) and translationally invariant (TI) partitions, (d) 
		mean-squared error decay as a function of compute time.}
		 \label{fig:results}
 \end{figure*}

We evaluate the effectiveness of the proposed approaches in
reconstructing a piecewise smooth function from noisy compressive
measurements.  In our simulations, the true signal (the black line in
Figs.~\ref{fig:results}(a--c)) is of length 1024.  We take 512 noisy
compressive measurements of the signal using a sensing matrix that
contains 32 randomly distributed nonzero elements per row.  This setup
yields a mean detector photon count of 50, ranging from as few as 22
photons, to as high as 94 photons.  We allowed each algorithm a fixed
time budget of three seconds in which to run, which is sufficient to
yield approximate convergence for all methods considered.  Each algorithm was
initialized at the same starting point, which was generated using a single E-step
of the EM algorithm.  That is if $z = A^Ty$, and $x : x_i = y_i/(Az)_i$, then we initialize with
$f^0 : f_i^0 = z_i (A^Tx)_i/(A^T \ones)_i$. The value
of the regularization parameter $\tau$ was tuned independently for
each algorithm to yield the minimal root-mean-squared error
$\text{RMS}=\|\hat{f} - f\|_2/\|f\|_2$ at the exhaustion of the
computation budget (see Table~\ref{tab:MSEs}).
\begin{table}[htbp]
   \centering
   \begin{tabular}{@{} ll @{}} 
      \toprule
      Reconstruction Method 	& RMS \\
      \midrule
      SPIRAL (TI) 		& 0.1427 \\
      SPIRAL (TV) 		& 0.1855 \\
      EM-MPLE (TI) 		& 0.2485 \\
      EM-MPLE (TV) 		& 0.2511 \\
      SPIRAL ($\ell_1$) 	& 0.2898 \\
      \bottomrule
   \end{tabular}
   \caption{Reconstruction accuracy as measured by 
   		the root-mean-squared error, $\text{RMS}=\|\hat{f} - f\|_2/\|f\|_2$.}
   \label{tab:MSEs}
   \vspace{-2ex}
\end{table}

\sloppypar{An examination of the results in Figs.~\ref{fig:results}(a-c)
  reveals that even though models within a partition (constant pieces)
  are less smooth than higher-order wavelets, this drawback is
  neutralized through a combination of cycle spinning (TI), the hereditary
  constraint, and a fast, non-iterative solver. 
We compare our algorithm with an EM algorithm based upon the same
maximum penalized likelihood estimation (MPLE) objective function
proposed in (\ref{eq:objective}), which has been used previously in
imaging contexts in which sparsity in a multiscale partition-based
framework was exploited \cite{willett:tmi03}.
Although
  Fig.~\ref{fig:results}(d) shows that the convergence behavior of the
  EM-MPLE approaches is
  more stable, their slow convergence ultimately hinders their
  performance as compared with the corresponding SPIRAL-based
  approaches. }
  In the case of SPIRAL-$\ell_1$, the estimate seems very
  oscillatory; a smoother estimate could be achieved by increasing the
  regularization parameter $\tau$, but this leads to an
  ``oversmoothed'' solution and increases the RMS of the estimate.
The large RMS values early in the SPIRAL iterations are due to
small values of $\alpha_k$ when defining $s^k$ in (\ref{eq:sk}),
which occur when the estimates are flat and
the consecutive iterates only change slightly.  
The Barzilai-Borwein method \cite{GPSR}, on which the choice of $\alpha_k$ is based,
is not monotone and can exhibit spikes in the iterates' RMS values 
(the non-monotone SpaRSA algorithm is also not immune to this behavior, for example).  
Although it is difficult to characterize the
convergence of this approach, it has been shown to be effective in solving
minimization problems (see \cite{GPSR} for details).
  The SPIRAL partition-based estimates appear to oscillate in
  steady state. This effect could be mitigated by forcing the
  objective function to decrease monotonically with iteration, as
  described in \cite{sparsa} and as in EM algorithms. However,
  empirically this appears to hinder performance, as convergence is
  significantly slowed.   We are currently investigating methods that use a
  non-monotonic approach to obtain a good approximation to the minimizer,
  followed by a monotonic method to enforce convergence.

\vspace{-.1cm}

\section{Conclusion}
\label{sec:conclusion}

We have developed computational approaches for signal reconstruction from photon-limited measurements---a situation prevalent in many practical settings.  Our method optimizes a regularized Poisson likelihood under nonnegativity constraints.  We have demonstrated that these methods prove effective in the compressed sensing context where an $\ell_1$ penalty is used to encourage sparsity of the resulting solution.  Our method improves upon current approaches in terms of reconstruction accuracy and computational efficiency.  By employing model-based estimates that utilize structure in the coefficients beyond that of a parsimonious representation, we are able to achieve greater accuracy with less computational burden.  Future work includes supplementing our algorithms with a debiasing stage, which maximizes the likelihood over the sparse support discovered in the main algorithm.  We also will be exploring the efficacy of alternative optimization approaches such as sequential quadratic programming and interior-point methods.

\bibliographystyle{IEEEbib}
\bibliography{strings,refs}

\end{document}